\documentclass{article}

\usepackage{amsthm, amssymb, amsmath, latexsym}
\usepackage[vcentermath,enableskew]{youngtab}

\newtheorem{prop}{Proposition}[section]
\newtheorem{thm}[prop]{Theorem}
\newtheorem{cor}[prop]{Corollary}
\newtheorem{lem}[prop]{Lemma}

\newcommand{\card}[1]{{\left\arrowvert #1 \right\arrowvert\,}}
\newcommand{\st}{{\,\arrowvert\, }}

\newcommand{\set}[1]{{\{ #1\}}}

\newcommand{\carp}[1]{{\mathcal #1 }}

\newcommand{\bo}{{\, \leq\, }}

\begin{document}

\title{Minimal and maximal elements in two-sided cells of $S_n$ and  Robinson-Schensted correspondence}

\author{{\bf Christophe Hohlweg}\\
Institut de Recherche Math\'ematique Avanc\'ee,\\
Universit\'e Louis Pasteur et CNRS,\\
7,rue Ren\'e Descartes\\
67084 Strasbourg, France\\
hohlweg@math.u-strasbg.fr}

\maketitle

\begin{abstract}
In symmetric groups, a two-sided cell is  the set of all
permutations which are mapped by the Robinson-Schensted
correspondence on a pair of tableaux of the same shape. In this
article, we show that the set of permutations in a two-sided cell
which have  a  minimal number of inversions   is the set of
permutations which have a maximal number of inversions in
conjugated Young  subgroups. We also give an interpretation of
these sets  with particular  tableaux, called reading tableaux. As
corollary,   we give  the set of elements  in a two-sided cell
which have a maximal number of inversions.
\end{abstract}

\noindent{\small{\it Author keywords:} Robinson-Schensted correspondence, number of inversions, two-sided cells .}

\section{Introduction}

In this article, we consider  the  symmetric group $S_n$. For
$w\in S_n$, {\it the length of} $w$, denoted  $\ell(w)$, is the
number of inversions of $w$.

The  Robinson-Schensted correspondence \cite{s} is the well-known
bijection $\pi : \, w \mapsto (P(w) , Q(w))$  between $S_n$ and
pairs of standard tableaux of the same shape (a partition of $n$).
For each partition $\lambda$ of $n$, we denote by $\carp
T^\lambda$ the set of all permutations which are mapped by $\pi$
on a pair of tableaux of shape $\lambda$. In the Kazhdan-Lusztig
theory, which we use in this article,  the sets $\carp T^\lambda$
are called  {\it two-sided cells} (see \cite{kaz,vogan,ariki}).

Our goal is to describe, for any partition $\lambda$ of $n$, the set
 $\carp{T}^{\lambda} _{min}$ of elements of minimal  length in
 $\carp T^\lambda$ and  the set $\carp T ^{\lambda} _{max}$ of elements
 of  maximal length in $\carp T^\lambda$.\\

For each composition $c=( n_1, \dots ,n_k)$ of $n$ (with $n_i \geq
1$), the {\it Young subgroup} $S_c = S_{n_1}\times \dots \times
S_{n_k}$ contains a unique permutation $\sigma_c$  of maximal
length. It is well-known that $\sigma_c$ is an involution, called
{\it the longest element} of $S_c$.  Denote $\lambda (c)$ the
decreasing reordering of $c$. It is well-known that two Young
subgroups $S_{c_1}$ and $S_{c_2}$ are conjugated in $S_n$ if and
only if $\lambda(c_1)=\lambda(c_2)$.

Sch\" utzenberger \cite{schutz2} has shown  that the map $T\mapsto
w_T= \pi^{-1}(T,T)$ is a bijection between the standard tableaux
of shape a partition of $n$ and the involutions of $S_n$ (see also
\cite{burge} where another interesting description of this
bijection is given).

A standard tableau $T$ is a \textit{reading  tableau} if it has
the following property: for any $1\leq p \leq n$, either $p$ is
in the first line of $T$, or if  $p$ is in the $i^{\textrm{th}}$
line of $T$ ($i>1$) then  $p-1$ is in the $(i-1)^{\textrm{th}}$
line of $T$.

As example, the \textit{column superstandard  tableaux} (which are
tableaux numbered from the bottom to the top of each column, from
left to right) are reading  tableaux. These particular tableaux
are the transposed of those defined by Garsia and Remmel in
\cite{garsiaremmel}.
 Our main result is the following:

\begin{thm}
\label{minimal}
Let $\lambda$ be a partition of $n$ and  $\carp{T}^{\lambda}$ be its associated two-sided cell, then
\begin{eqnarray*}
\carp T ^{\lambda} _{min}&=&\set{\sigma_c \st \lambda(c)= \lambda^{t}}\\
&=& \set{w_T \st T\textrm{ is a reading  tableau of shape
}\lambda},
\end{eqnarray*}
where $\lambda^t$ denotes the {\it conjugate partition} of
$\lambda$.
\end{thm}

\noindent\textit{Example.} Consider the partition $\lambda = (3 ,2
, 1 , 1)$ of $7$; $\lambda^t = (4 , 2 , 1)$. Then the reading
 tableaux of shape $\lambda$  are
$$
T_1 = \young(157,26,3,4)\quad\equiv\quad
\young(4::,3::,2::,1::,:6:,:5:,::7) %
\quad ; \quad %
T_2 =
\young(137,24,5,6)\quad\equiv\quad\young(2::,1::,:6:,:5:,:4:,:3:,::7)
$$
$$
T_3 =
\young(126,37,4,5)\quad\equiv\quad\young(1::,:5:,:4:,:3:,:2:,::7,::6)%
\quad ; \quad%
T_4 =
\young(156,27,3,4)\quad\equiv\quad\young(4::,3::,2::,1::,:5:,::7,::6)
$$
$$
T_5 = \young(134,25,6,7)\quad\equiv\quad\young(2::,1::,:3:,::7,::6,::5,::4)%
\quad ; \quad%
T_6 =
\young(124,35,6,7)\quad\equiv\quad\young(1::,:3:,:2:,::7,::6,::5,::4)
$$
The right skew tableau is taken in the plactic class of the
corresponding tableau viewed in the plactic monoid
\cite{lascoux-schutz} (see also \cite{fulton}). The corresponding
involutions are then
$$
w_{T_1}= 4\,3\, 2\,1\,6\, 5\, 7 = \sigma_{(4,2,1)}\quad ; \quad w_{T_2}= 2\,1\,6\, 5\,4\,3\,7 = \sigma_{(2,4,1)}
$$
$$
w_{T_3}= 1\, 5\,4\,3\,2 \,7 \,6= \sigma_{(1,4,2)} \quad ; \quad w_{T_4}= 4\,3\,2\,1\, 5\, 7\,6  = \sigma_{(4,1,2)}
$$
$$
w_{T_5}= 2\,1\, 3\,7\, 6\, 5\,4 = \sigma_{(2,1,4)} \quad ; \quad w_{T_6}= 1\, 3\,2 \,7\,6\, 5\, 4 = \sigma_{(1,2,4)}
$$
and $\carp T^{(3,2,1,1)}_{min} = \set{ \sigma_{(4,2,1)}, \sigma_{(2,4,1)}, \sigma_{(1,4,2)},
\sigma_{(4,1,2)}, \sigma_{(2,1,4)} , \sigma_{(1,2,4)} }$.\\
\\

Let $\lambda$ be a partition of $n$ and  $T$ be a standard tableau
of shape $\lambda$. The \textit{Sch\"utzenberger evacuation} of
$T$, denoted by $ev \, (T)$, is a tableau of shape $\lambda$
\cite{sc} (see also \cite[p.128-130]{sagan}). The evacuation
illustrates the conjugation and the left (and right)
multiplication by the longest element $\sigma_{(n)}$ in $S_n$.  In
particular  $Q (w \sigma_{(n)}) = ev\,( Q(w)^t)$, for any $w\in
S_n$, and $\carp T ^{\lambda} \sigma_{(n)} =\sigma_{(n)}\carp T
^{\lambda}= \carp T ^{\lambda^t}$. Denote
$d_c=\sigma_{(n)}\sigma_c$. As $\ell(\sigma_{(n)}w) =
\ell(\sigma_{(n)})-\ell(w)$, we obtain the following corollary:

\begin{cor}\label{maximum}  Let $\lambda$ be a partition of $n$ and  $\carp{T}^{\lambda}$
be its associated two-sided cell,  then
\begin{eqnarray*}
\carp T ^{\lambda} _{max}&=&\set{d_c \st \lambda(c)= \lambda}\\
&=& \set{ w \st ev\,( Q(w)^t)= P(w\sigma_{(n)})\textrm{ is a
reading  tableau of shape }\lambda^t}.
\end{eqnarray*}
\end{cor}

In the theory of Coxeter groups, the element $d_c$ is well-known
as the unique element of maximal length in the set of minimal
right coset representatives of $S_c$ \cite[Chapter 2]{geck}.

As a by-product of our proof,  we obtain,  in Section
\ref{section-duflo}, that if $w$ is an involution,  the
Kazhdan-Lusztig polynomial $P_{e,w}=1$ if and only if $w$ is the
longest element of a Young subgroup, where $e$ denotes the
identity of $S_n$.  More precisely, we show that an involution $w$
avoids the pattern $3412$ and $4231$ if and only if $w$ is the
longest
element of a Young subgroup of $S_n$. \\

To our knowledge, our  results are the first results relating the
Robinson-Schensted transformation and the length function (number
of inversions) of a permutation. There is no evident link between
both. Our proof is non combinatorial and uses heavily the
$a-$function of Lusztig \cite{lusztig,lusztig-cells} (questions
about the leading term of Kazhdan-Lusztig polynomials and the
$a-$function are heavily studied, see for instance
\cite{xi1,xi2}). Trying to find a combinatorial proof (which is a
challenge) leads first to the following difficulty: for any
permutations $w, x$  in a two sided-cell $\carp T^\lambda$, there
are permutations $w_1, \dots, w_k \in \carp T^\lambda$ such that
$w_1 = w$, $w_k = x$ and $w_{i+1}$ is obtained from $w_{i}$ by a
Knuth or a dual-Knuth relation (see \cite{sagan,fulton}). But  a
permutation $w$ may be `locally minimal', that is whenever $w$ (or
$w^{-1}$) admits a Knuth or a dual-Knuth elementary relation, this
relation increases the length. The involution $\sigma=632541 \in
S_6$ is locally minimal, but not of minimal length, in its
two-sided cell.

It would be interesting to find a  purely combinatorial proof of
the main result. It is apparently an open problem to read the
length of a permutation $w$ directly on the pair of tableaux $\pi
(w)$ (however, see \cite{signature}, where the author gives a way
to read the signature on the pair of tableaux). Fortunately, the
Lusztig $a$-function gives us a way to avoid this problem.

%
%
%
%
%
%
%
%

\section{Consequences of the main result}\label{kl-cells}

We denote  a partition of $n$ by $\lambda = (\lambda_1 , \dots ,
\lambda_k)$, with $\lambda_1 \geq \dots \geq \lambda_k \geq 1$.
Our reference for the general theory of the symmetric group is
\cite{sagan}.

For any partition $\lambda = (\lambda_1 , \dots , \lambda_k )$ of
$n$, we  define, for $i> 0$,
$$
m_i (\lambda) = \card{\set{j \st \lambda_j = i}} .
$$
The number $m_i(\lambda)$ is called the \textit{multiplicity} of
$i$ in $\lambda$ (see \cite{macdo}). Observe that $m_i(\lambda)=0$
for all $i>n$, since  $\sum_{i} \lambda_i =n$. It is well-known
that the multinomial coefficient
$$
\left(
\begin{array}{c}
m_1(\lambda)+m_2(\lambda)+\dots +m_n(\lambda)\\
m_1(\lambda),m_2(\lambda), \dots ,m_n(\lambda)
\end{array} \right)
$$
is the number of compositions associated to $\lambda$. Hence, we obtain the following corollary.

\begin{cor} Let $\lambda$ be a partition of $n$, then
$$
\card{\carp T ^{\lambda} _{min}} = \left(
\begin{array}{c}
m_1(\lambda^t)+m_2(\lambda^t)+\dots +m_n(\lambda^t)\\
m_1(\lambda^t),m_2(\lambda^t), \dots ,m_n(\lambda^t)
\end{array} \right) ,
$$
which is the number of compositions $c$ of $n$ such that $\lambda(c)=\lambda^t$.
\end{cor}

 The minimal elements in two-sided  cells are linked to another important number in combinatorics
$$
n(\lambda) = \sum_{i=1}^k
\left(\begin{array}{@{}c@{}}
\lambda_i^t\\
2
\end{array}\right) ,
$$
see \cite[p.2-3]{macdo}.

\begin{cor}
\label{length_minimal}
Let $\lambda$ be a partition of $n$ and  write $\lambda^t =
(\lambda_1^t , \dots , \lambda_k^t)$. Then $\ell(w) = n(\lambda)$
for all $w\in\carp T ^{\lambda} _{min}$.
\end{cor}

\begin{proof}
Let $c$ be a composition of $n$ such that $\lambda (c)
=\lambda^t$. Then $\ell(\sigma_c) =\ell(\sigma_{\lambda^t}$. Let
$w_{i}$ be the longest element of the Young subgroup
$S_{\lambda^t_i}$, then  $\ell(w_{i}) =
\left(\begin{array}{@{}c@{}}
\lambda^t_i\\
2
\end{array}\right)$. Therefore
$$
\ell(\sigma_{\lambda^t}) = \sum_{i=1}^k
\left(\begin{array}{@{}c@{}}
\lambda^t_i\\
2
\end{array}\right)
$$
since $\sigma_{\lambda^t}= w_{1}\dots w_{k}$ (seen as a word on
the letters $1,\dots , n$) and that the letters in $w_{i+1}$ are
greater than the letters in $w_{i}$. The corollary follows from
Theorem~\ref{minimal}.
\end{proof}

As in the case of minimal elements, we have the following
corollaries:

\begin{cor} Let $\lambda$ be a partition of $n$, then
$$
\card{\carp T ^{\lambda} _{max}} = \left(
\begin{array}{c}
m_1(\lambda)+m_2(\lambda)+\dots +m_n(\lambda)\\
m_1(\lambda),m_2(\lambda), \dots ,m_n(\lambda)
\end{array} \right) ,
$$
which is the number of compositions $c$ of $n$ such that $\lambda(c)=\lambda$.
\end{cor}

\begin{cor} Let $\lambda= (\lambda_1 , \dots , \lambda_k)$ be a partition of $n$,  then
$$
\ell(w) =\left(\begin{array}{@{}c@{}}
n\\
2
\end{array}\right) - \sum_{i=1}^k
\left(\begin{array}{@{}c@{}}
\lambda_i\\
2
\end{array}\right)
$$
for all $w\in\carp T ^{\lambda} _{max}$.
\end{cor}

\begin{proof} As $\ell(\sigma_{(n)}) =\left(\begin{array}{@{}c@{}}
n\\
2
\end{array}\right)$ and
$\ell(\sigma_{(n)}  w ) = \ell(\sigma_{(n)}) - \ell(w)$, for any
$w\in S_n$, the corollary follows from  same argument than in the
proof of Corollary~\ref{length_minimal}.
\end{proof}

%
%
%
%
%
%
%
%
%
%

\section{Proof of Theorem~\ref{minimal}}\label{def_minimal}\label{s2} \label{rc-tableau}\label{section-duflo}

The following proposition implies that
$$
\set{\sigma_c \st \lambda(c)= \lambda^{t}}= \set{w_T \st T\textrm{
is a reading  tableau of shape }\lambda}\subset\carp T ^{\lambda}
.
$$

\begin{prop}\label{cctapis1} Let $\lambda$ be a partition of $n$; then the following conditions are equivalent:
\begin{enumerate}
\item[i)] $T$ is a reading  tableau of shape $\lambda$; \item[ii)]
$w_T=\sigma_c$, where $c$ is a composition of $n$ such that
$\lambda(c)=\lambda^{t}$.
\end{enumerate}
\end{prop}

\begin{proof}
Recall that the longest element of a Young subgroup is an involution, since it is unique.

Assume $(i)$. As $T$ is a reading tableau, if $n$ is in the row
$T_i$, one has  $1\leq p \leq n-1$ such that $p+1$ is in the first
row of $T$,  $p+i=n$ and $p+j$ is at the end of the row $T_{j}$,
for all $1\leq j \leq i$. One applies the $i$ first steps of the
inverse of  Robinson-Schensted correspondence, hence
$$
w_T = w_{T'}\ n \ \dots \ p+1 .
$$
where $T'$ is the standard Young tableau obtained by deleting
$p+1, \dots ,n$ in $T$. Thus $w_{T'}$ is a permutation on the set
$\set{1,\dots , p}$. Observe that $T'$ is also a reading tableau.
The shape of $T'$ is denoted by $\lambda'$. By induction on $n$,
$w_{T'}$ is the longest element of the Young  subgroup $S_{c'}$,
where $\lambda(c') = \lambda'$. Then $w_T$ is the longest element
of the Young  subgroup $S_{c'}\times S_{i}$. Let $c =(c' , i)$; it
is now easy to see that  $\lambda(c)=\lambda^t$.

Conversely, let $c = (n_1 ,\dots n_k)$ and use induction and
similar arguments with direct Robinson-Schensted correspondence on
the permutation $$w_T = n_1 \ \dots\ 1\  w' ,$$ where $n_1 \
\dots\ 1$ is the longest element of the Young subgroup $S_{n_1}$
and $w'$ is the longest element of the Young subgroup
$S_{n_2}\times \dots\times S_{n_k}$.
\end{proof}

Now, it remains to prove that $\set{\sigma_c \st \lambda(c)=
\lambda^{t}}=\carp T^\lambda_{min}$, to end the proof of
Theorem~\ref{minimal}.

\paragraph{The Lusztig $a$-function:} We consider  the  symmetric group $S_n$ as a Coxeter system
$(W,S)$ of type $A_{n-1}$ with $W=S_n$ and  generating set $S$
consisting of the $n-1$ simple transpositions $\tau_i=(i,\ i+1)$,
where  $i=1, \dots , n-1$. Then $\ell(w)$ is also the length of
$w$ as a reduced word in the elements of $S$. A classical
bijection between subsets of $S$ and compositions of $n$ is
obtained as follow: Let $I\subset S$ and $S\setminus I =
\set{\tau_{i_1}, \dots ,\tau_{i_k}}$ with $1\leq i_{1} <i_2 <\dots
<i_k \leq n-1$. Set $n_1 = i_1$, $n_2 = i_2 - i_1 +1$, $\dots$,
$n_k = n - i_k$, then $n_i$ are non-negative integers. By this
way, we have obtained a unique composition  $c_I = (n_1 , \dots ,
n_k)$ of $n$ associated to $I$. Moreover,
$$
W_{I} = S_{n_1}\times \dots \times S_{n_k}.
$$
Therefore, as is well-known the Young subgroups of $S_n$ are
precisely the \textit{parabolic subgroups} of $S_n$ (see
\cite[Proposition 2.3.8]{geck}).

Our basic references for the work of Kazhdan and Lusztig are
\cite{kaz}, \cite{lusztig} (see also \cite{curtis}).  We denote by
$\bo$ the Bruhat order on $S_n$.

Let $\mathcal A = \mathbb Z [q^{1/2},q^{-1/2} ]$ where $q^{1/ 2}$
is an indeterminate. Let $\mathcal H$  be the Hecke algebra over
$\mathcal A$ corresponding to $S_n$.  Let $(T_w)_{w\in S_n}$ be
the standard basis of $\mathcal H$ and $(\widetilde{T}_w)_{w\in
S_n}$ the basis defined as follows:
$$
\widetilde{T}_w=q^{-\ell(w)/ 2}T_w .
$$

In \cite[Theorem 1.1]{kaz}, Kazhdan and Lusztig have shown that
there is a basis $(b_w)_{w\in S_n}$ of $\mathcal H$, called the
\textit{Kazhdan-Lusztig basis},  such that
$$
b_w = \sum_{y\leq w} (-1)^{\ell(w) - \ell(y)} q^{(\ell(w) - \ell(y))/2} P_{y,w}(q^{-1})\widetilde{T}_y ,
$$
where $P_{y,w}\in \mathcal A$ are  the Kazhdan-Lusztig
polynomials. Moreover, they have defined three   equivalence
relations on $S_n$, with equivalence classes that are called
\textit{left cells}, \textit{right cells} and \textit{two-sided
cells}. In our case,   the following result of Vogan and Jantzen
result on $S_n$ \cite{jantzen,vogan} gives the link with the
Robinson-Schensted correspondence (see also \cite{ariki}): the set
$\carp T ^{\lambda}$ is a two-sided cell for all partitions
$\lambda$ of $n$; and any two-sided cell of $S_n$ arises by this
way.

Following Lusztig  \cite{lusztig-cells, lusztig}, let $h_{x,y,w}$
be the  structure constants of the Kazhdan-Lusztig base
$(b_w)_{w\in W}$, that is
$$
b_x b_y = \sum_{w\in W} h_{x,y,w}\, b_w .
$$

Denote $\delta(w)$ the degree of the Kazhdan-Lusztig polynomial
$P_{e,w}$ as a polynomial in $q$. Write $u=q^{1/ 2}$. Let $a(w)$
be the smallest integer such that for any $x,y\in S_n$,
$u^{a(w)}h_{x,y,w}\in \mathcal{A}^{+}$, where
$\mathcal{A}^{+}=\mathbb{Z}[u]$ (this is well defined for any Weyl
group). In \cite{lusztig-cells, lusztig}, Lusztig has shown the
following properties about the $a-$function:
\begin{enumerate}
\item[a)] $a(w)\leq \ell(w) -2 \delta(w)$ (\cite[Section 1.3]{lusztig});

\item[b)] The $a-$function is constant on two-sided cells (\cite[Theorem 5.4]{lusztig-cells}).

\item[c)] For any $I\subset S$, $a(\sigma_{c_I})=\ell(\sigma_{c_I})$
(\cite[Corollary 1.9 (d) and Theorem 1.10]{lusztig}).
In other words, for any  composition $c$ of $n$, $a(\sigma_c)=\ell(\sigma_c)$.

\item[d)] Let $\mathcal{D}=\set{w\in W \st a(w)= \ell(w) -2
\delta(w) }$, then each element in $\mathcal D$ is an involution,
called a \textit{Duflo involution} (\cite[Proposition
1.4]{lusztig}). In symmetric groups, all involutions are Duflo
involutions. Indeed, each left cell contains a unique Duflo
involution \cite{lusztig};   left cells are precisely coplactic
classes (the sets of permutations having the same right tableau
under $\pi$, see for instance \cite{ariki}), and each coplactic
class contains a unique involution.
\end{enumerate}

Let $\lambda$ be a partition of $n$ and $\carp T^\lambda$ be its
associated two-sided cells. Properties (b) and (c) imply that
$a_\lambda := a(\sigma_{\lambda^t}) = a(w)$, for all $w\in \carp
T^\lambda$. Therefore, by (a),
$$
\ell(\sigma_c) =a_\lambda =a(w) \leq \ell(w) ,
$$
for any $w\in \carp T^\lambda$. Thus
$$
\set{\sigma_c \st \lambda(c)= \lambda^{t}}\subset \carp T ^{\lambda} _{min} .
$$

Now, let $w\in \carp T ^{\lambda} _{min}$, then $a(w)=a_\lambda =
\ell(\sigma_{\lambda^t}) = \ell(w)$, since $\sigma_{\lambda^t}\in
\carp T ^{\lambda} _{min}$. Property (d) implies that $w$ is a
Duflo involution and $\delta (w) =0$.

By Proposition~\ref{cctapis1}, $\sigma_c \in \carp T^\lambda$
implies $\lambda (c) = \lambda^t$ Therefore, Theorem~\ref{minimal}
is a direct consequence of the  following result,  which  gives a
surprising criterion about  the degree $\delta(w)$  of the
Kazhdan-Lusztig polynomial $P_{e,w}$, for $w\in S_n$ an
involution.

\begin{prop}\label{min1} Let $w\in S_n$, then the following conditions are equivalent:
\begin{enumerate}
\item[i)] $w$ is an involution and  $\delta(w)=0$;
\item[ii)] $w=\sigma_c$, for some composition $c$ of $n$.
\end{enumerate}
\end{prop}

\paragraph{KL Polynomials and smoothness of Schubert Varieties:} We say that a permutation
$w\in S_n$, seen as a word $w=x_1 \dots
x_n$, \textit{avoids the pattern} $4231$ (resp. \textit{avoids the
pattern} $3412$)  if there is no $1\leq i<j<k<l\leq n$ such that
$x_l < x_j < x_k <x_i$ (resp. $x_k <x_l <x_i <x_j$). In other
words, there is no subword of $w$ with the same relative order as
the word $4231$ (resp. $3412$).

Here, we  link these definitions with Kazhdan-Lusztig polynomials
by the way of the following well-known criterion: Let  $w\in S_n$,
then
$$
(\diamond) \qquad P_{e,w}=1 \iff w \textrm{ avoids the patterns } 4231 \textrm{ and }3412.
$$

Indeed, on one hand,  Lakshmibai and Sandhya have shown that a
Schubert variety $X(w)$, $w\in S_n$, is smooth if and only if $w$
avoids the pattern $3412$ and $4231$ (\cite{lake-sand} or see
\cite[Theorem 8.1.1]{billey}).

On the other hand, Deodhar \cite{deodhar85} has shown a useful
characterization of the smoothness by the way of Kazhdan-Lusztig
polynomials: Let $w\in S_n$ then $P_{e,w}=1$ if and only if $X(w)$
is smooth.

\paragraph{Proof of Proposition~\ref{min1}:} By the above  discussion,
Proposition~\ref{min1} is a direct consequence of the following lemma.

\begin{lem}\label{schubert} Let $w\in S_n$ an involution, then the following statements are equivalent
\begin{enumerate}
\item[i)] $w$ avoids the patterns $4231$ and $3412$;
\item[ii)] there is a composition $c$ of $n$   such that $w = \sigma_{c}$;
\end{enumerate}
\end{lem}
\begin{proof} $(ii)\Rightarrow (i)$: write $c=(c_1,\dots,c_k)$. If
$k=1$ then $\sigma_c=\sigma_{(n)}$ avoids the patterns $4231$ and
$3412$. If $k>1$ then $\sigma_c$ is the image of
$(\sigma_{(c_1)},\dots,\sigma_{(c_k)})$ under the canonical
isomorphism between $S_{c_1}\times\dots\times S_{c_k}$ and $S_c$.
Conclude by induction on $n$.

$(i)\Rightarrow (ii)$: one sees $w=x_1 \dots x_n$ as a word on the
letters $1,\dots , n$.

One proceeds by induction on $n$. Therefore, one may suppose that
$(i)\Rightarrow (ii)$ for all proper Young subgroups of $S_n$. If
$n \leq 4$, it is readily seen. Suppose $n>4$.

If $x_1 = 1$, then $w\in S_1 \times S_{n-1}$, and the lemma
follows by induction.

If $n>x_1 = p >1$, then $x_p = 1$ and $1\leq x_i \leq p$, for all
$1\leq i\leq p$. Otherwise, there is $1<i<p$ such that $x_i >p$.
In other words, there is $1<i<p<x_i$ such that $x_p = 1 <
x_{x_i}=i <x_1 = p < x_i$, that is, $w$ has  the pattern $3412$
which is a contradiction.

Hence $w\in S_p \times S_{n-p}$ and the lemma follows by induction.

If $x_1=n$, then $x_n = 1$ one just has to show that $w=w_0$.
Otherwise, there is $1< i< n-1$ such that $x_i < x_{i+1}$  (since
if $i=1$,  $x_1 =n <x_2$ and if $i+1=n$, $x_{n-1}<x_n=1$ which are
contradictions). Thus there is $1<i<i+1<n$ such that $x_n<x_i
<x_{i+1}<x_1$, that is, $w$ has  the pattern $4231$ which is a
contradiction.
\end{proof}

\subsection*{Acknowledgements}

I am  indebted to C.~Chauve and C.~Reutenauer  for useful
discussions about Robinson-Schensted correspondence and plactic
classes  and to P.~Polo for his useful remarks. I would thank
P.~Baumann, R.~B\'edard and F.~Chapoton  for pointing out some
references.

\end{document}